\documentclass[11pt]{article} \usepackage{amsfonts}  

\makeatletter
\def\quality{\textheight=240mm \textwidth=160mm \topmargin=0Truein
             \ifcase \@ptsize \hoffset=-23mm
                     \or \hoffset=-20mm \or \hoffset=-15mm \fi}
\def\bdraft{\pagestyle{myheadings} 
           \textheight=10.5truein \textwidth=7.5truein \parindent=8pt
           \voffset=-1truein \topmargin=30pt \headheight=10pt \headsep=3pt
           \ifcase \@ptsize \hoffset=-1.5truein \or \hoffset=-1.35truein
                        \or \hoffset=-1.15truein \fi}
\makeatother
\quality

\def\cyl#1{\langle#1\rangle}
  \def\cM{{\cal M}}

\def\IR{{\mathbb{R}}} \def\IZ{{\mathbb{Z}}} \def\t{\tilde}
\def\cB{{\cal B}}  
\def\cM{{\cal M}}   \def\cA{{\cal A}}
\def\toas#1{\stackrel{#1}{\longrightarrow}}

\def\cyl#1{\langle#1\rangle}              \def\CR{$$ $$}
\def\E{\mathbb{E}} \def\P{\mathbb{P}}  
\def\dist{\varrho} \def\den{\rho}  \def\comment#1{}

\def\bline(#1,#2)(#3,#4)(#5){\put(#1,#2){\line(#3,#4){#5}}}  

\newcommand\mlbscale{1pt} 
\newif\iffigs\figstrue 

\def\bfig(#1,#2)#3#4{\begin{figure} \begin{center}
    \framebox{\setlength{\unitlength}{\mlbscale}
       \iffigs \begin{picture}(#1,#2) #3 \end{picture}
       \else \begin{picture}(60,10)(0,0)
                   \put(0,0){\framebox(60,10){Figure}} \end{picture} \fi}
    \end{center} \caption{#4} \end{figure}}

\def\Bfig(#1,#2)#3#4{\begin{figure} \begin{center}
    \setlength{\unitlength}{\mlbscale}
       \iffigs \begin{picture}(#1,#2) #3 \end{picture}
       \else \begin{picture}(60,10)(0,0)
                   \put(0,0){\framebox(60,10){Figure}} \end{picture} \fi
    \end{center} \caption{#4} \end{figure}}

 \def\IZ{{\mathbb{Z}}}   \def\ep{\varepsilon}

\def\function#1{\left\{\!\!\!\begin{array}{ll} #1 \end{array} \right.}
\def\bea#1{\begin{eqnarray*} #1 \end{eqnarray*}} \def\a{\!\!\!&\!\!\!\!&}
\def\beq#1#2{\begin{equation} \label{#1} #2 \end{equation}}

\def\thname{Theorem}     \def\lmname{Lemma}      \def\prname{Proposition}
\def\dfname{Definition}  \def\crname{Corollary}  \def\rmname{Remark}

\newtheorem{theorem}{\thname}[section]   
\newtheorem{lemma}{\lmname}[section]     
\newtheorem{corollary}[lemma]{\crname}   

\newtheorem{dftn}{\dfname}[section]
\newtheorem{rmrk}[lemma]{\rmname}

             \catcode`@=11 \@addtoreset{equation}{section} \catcode`@=12

\def\proof{\smallskip \noindent {\bf Proof. \ }}
\newcommand\emptysquare{\ \ $\Box$}
\def\qed{\hfill\emptysquare\linebreak\smallskip\par}

\begin{document}

\title{
       On quasi successful couplings of Markov processes}
\author{\large Michael Blank\thanks{Russian Academy of Sci.,
            Inst. for Information Transm. Problems,
            and Observatoire de la Cote d'Azur,
            ~e-mail: blank@iitp.ru}, ~~
        Sergey Pirogov\thanks{Russian Academy of Sci.,
            Inst. for Information Transm. Problems,
            ~e-mail: pirogov@mail.ru}}
\date{October 3, 2006} 
\maketitle

\maketitle

\begin{abstract} The notion of a successful coupling of Markov
processes, based on the idea that both components of the coupled
system ``intersect'' in finite time with probability one, is
extended to cover situations when the coupling is unnecessarily
Markovian and its components are only converging (in a certain
sense) to each other with time. Under these assumptions the unique
ergodicity of the original Markov process is proven. A price for
this generalization is the weak convergence to the unique
invariant measure instead of the strong one. Applying these ideas
to infinite interacting particle systems we consider even more
involved situations when the unique ergodicity can be proven only
for a restriction of the original system to a certain class of
initial distributions (e.g. translational invariant ones).
Questions about the existence of invariant measures with a given
particle density are discussed as well.
\end{abstract}

\section{Introduction}

Let $(X,\cB)$ be a measurable space equipped with a metric
$\dist(\cdot,\cdot)$ and assume that the resulting metric space is
locally compact. A Markov chain $\xi^t$ acting on this phase space
is defined by a family of transition probabilities $P^t(x,A)$ to
jump from a point $x\in X$ to a measurable set $A\in\cB$ during
the time $t\ge0$ and an initial distribution of the random
variable $\xi^0$ representing the initial state of the Markov
chain. The Markov chain $\xi^t$ generates two semigroups of
operators $P^t\phi(x):=\int\phi(y)P^t(x,dy)$ acting on bounded
measurable functions and $P_*^t\mu(\phi):=\mu(P^t\phi)$ acting on
measures.

One of the first questions in the analysis of Markov chains is the
existense/uniqueness of their invariant (stationary) measures
(solutions to the equation $P_*^t\mu=\mu~~\forall{t}$) and
convergence to them for various initial distributions. We shall
describe a novel approach to study these questions applicable to a
reasonably broad class of Markov chains.

To a large extent various coupling results about Markov chains are
based on the so called {\em coupling inequality} applied as
follows: let $\xi^t$ and $\acute\xi^t$ be two Markov chains with
the same transition probabilities and initial conditions $\xi^0=x,
~\acute\xi^0=y$ defined on a common probability space
$(\Omega,{\cal F},\P_{x,y})$ where the distribution $\P_{x,y}$
depends on $x,y$ measurably. Denote by $\tau^0$ a random
variable, called {\em intersection time}, equal to the first
moment of time such that $\dist(\xi^t,\acute\xi^t)=0 ~~\forall
t\ge\tau^0$. In other words $\tau^0$ is the moment of
intersection of realizations of two Markov processes started from
the points $x,y\in X$. Using this notation the coupling
inequality can be written as follows:
\beq{e:coupling-ineq}
    {\sup_{A\in\cB}|\P_{x,y}(\xi^t\in A)-\P_{x,y}(\acute\xi^t\in{A})|
     \le \P(\tau^0>t) .}%

Assuming now that there exists a probability {\em invariant
measure} $\mu$ (called also a {\em stationary distribution}) of
the Markov chain $\xi^t$ and that the process $\acute\xi^t$ is
started in the stationary distribution $\mu$, we have that the
probability that $\acute\xi^{t}\in A$ is equal to $\mu(A)$ for
any set $A\in\cB$ and $t\ge0$.

Thus if the intersection time is finite with probability one the
inequality~(\ref{e:coupling-ineq}) proves the convergence in total
variation of the law of $\xi^t$ to $\mu$ (see
Section~\ref{s:functional}). We refer the reader for excellent
reviews of various results based on the coupling inequality to
\cite{Lig,Num,Thor}.

Recall that a {\em coupling}\/ is an arrangement of a pair of
processes on a common probability space to facilitate their direct
comparison, namely the coupling of Markov chains $\xi^t$ and
$\acute\xi^t$ is a pairs process $(\xi^t,\acute\xi^t)$ defined on
the direct product space $X\times{X}$ and satisfying the
assumptions %
$$ \P_{x,y}((\xi^t,\acute\xi^t)\in A\times X)
     = \P_{x}(\xi^{t}\in A), \CR
   \P_{x,y}((\xi^t,\acute\xi^t)\in X\times A)
     = \P_{y}(\acute\xi^{t}\in A) $$
for each $t\ge0, ~A\in\cB$. Extending this definition for arbitrary 
initial distributions $\mu,\nu$ instead of initial points $x,y$ one 
gets the coupling in this case as well. 

One says that a coupling is {\em successful}\/ if the components 
of the coupled processes coincide starting from the finite moment 
of the intersection time. Under certain assumptions on the Markov 
chain $\xi^t$ it is believed that the existence of the successful 
coupling is equivalent to the convergence of distributions 
$\P(\xi^t\in\cdot)$ and $\P(\acute\xi^t\in\cdot)$ 
(see, e.g. \cite{Num}). Of course, such equivalence cannot hold 
for arbitrary Markov chains. Let $\xi^t$ is described by a 
deterministic difference relation $x\to{x/2}$.
This relation defines a one-dimensional Markov chain with the
unique invariant measure $\mu$ -- the delta-measure at the
origin. Moreover, any probability measure converges weakly under
the action of this Markov chain to $\mu$. On the other hand, two
realizations of this Markov chain started at nonzero points $x,y$
will never intersect in finite time and thus no successful
coupling is possible.\footnote{One might argue that in this toy
model all ``states'' of the phase space except for the origin are
``unessential'' and this is the main reason of such behavior. This
is not the case and we shall return to this question in the
analysis of infinite particle systems in
Section~\ref{s:particles}.} %
Nevertheless any two realizations are becoming arbitrary close to
each other, which leads to a generalization of the notion of the
successful coupling.

Fix a certain coupling of the Markov chains $\xi^t$ and
$\acute\xi^t$ and for a given $\ep>0$ denote by $\tau^\ep$ a
random variable, called {\em quasi intersection time}, equal to
the first moment of time such that %
$\dist(\xi^t,\acute\xi^t)\le\ep ~\forall t\ge\tau^\ep$.

We shall say that the coupling is {\em metrically quasi successful}\/ if %
\beq{e:m-q-success}{\P_{x,y}(\tau^\ep<\infty)=1 
                                       \quad\forall x,y\in X, ~\ep>0;}%
and {\em topologically quasi successful}\/ if %
\beq{e:t-q-success}{
     \P_{x,y}(1_A(\xi^t)\ne1_A(\acute\xi^t))
     \toas{t\to\infty}0  \quad\forall x,y\in X} %
and any open set $A\in\cB$. Here $1_A(\cdot)$ stands for the
indicator function of the set $A\in\cB$.

These definitions represent quite natural generalizations of the
usual successful coupling. To the best of our knowledge the only
published result in this direction is the discussion of the so
called shift-coupling with close enough intervals of time in
\cite{Thor}, see also some particular results about Markovian
couplings for TASEP in \cite{EG,Lig}.

Both definitions of the quasi successful coupling look very
similar, but neither of them implies another one. Indeed, in the
above example of the coupling of two copies of the processes
defined by the map $x\to x/2$ the coupling is clearly metrically
quasi successful. To show that it is not topologically quasi
successful set $A:=\IR\setminus\{0\}, x=0,y=1$. Then $1_A(\xi^t)\equiv0$ 
while $1_A(\acute\xi^t)=1~\forall t\ge0$, which implies %
$\P_{x,y}(1_A(\xi^t)\ne1_A(\acute\xi^t))=1$. On the other hand,
it is obvious that the property (\ref{e:t-q-success}) cannot
imply (\ref{e:m-q-success}) since arbitrary large deviations
occurring with small probabilities are allowed in
(\ref{e:t-q-success}) but not in (\ref{e:m-q-success}).

The coupling we consider needs not to be Markovian, however in the
latter case the arguments might be significantly simplified (see
Section~\ref{s:markov}). Another important question also
discussed in Section~\ref{s:markov} is whether it suffices to
have two coupled Markov chains become equal (close enough) at a
single moment of time, or whether it is necessary to have them
remain equal for all future times as it is assumed in the
definition of the (quasi) intersection time.

The situation is becoming even more involved when we control the
convergence between the components of a coupled system only up to
a semimetric rather than a metric. This is typical for Markov
chains described infinite interacting particle systems (see
Section~\ref{s:particles} and \cite{Lig} for a general review). In
such cases there is no hope to study all invariant measures but as
we shall show, using the ideas developed in the analysis of quasi
successful couplings, one can prove the uniqueness of an invariant
measure for a system restricted to a reasonably large class of
initial distributions (e.g. translational invariant ones).

It is worth note that various models of {\em continuous time}
infinite interacting particle systems are very thoroughly studied
(see, e.g. \cite{Lig} and further references therein). Methods
used in these references are very model specific and are based on
rather different arguments (monotonicity, positivity, duality,
etc.) instead of a direct application of a version of a successful
coupling. Moreover, there are only a few mathematical results
\cite{GG, ERS, Grif, DKT, Bl-erg, Bl-hys, BF}\footnote{
   Compare this to a rather extensive list of publications
   related to the continuous time case in \cite{Lig}.
   Note also that the last three  references are dedicated to
   a pure deterministic setting.} %
about {\em discrete time} infinite interacting particle systems.
The main reason here is that in the continuous time case with
probability one only one ``interaction'' between particles can
happen at a given time, while in the discrete time setting an
arbitrary (even infinite) number of ``interactions'' may happen
simultaneously.

One of the main purpose of the present paper is to fill this gap
and to develop some direct tools to study the latter situation.
Therefore we shall discuss mainly the discrete time case,
explaining (where it is necessary) how to change the arguments to
prove the corresponding continuous time versions of our arguments.
The construction of couplings satisfying our assumptions for
specific particle systems is not trivial and will be discussed in
a separate paper.

Throughout this paper we shall not discuss the existence of
invariant measures for Markov processes under consideration
except for a special case of invariant measures with a given
particle density discussed in Section~\ref{s:inv-mes-density}.
In general to prove the existence one might need to assume
that these processes satisfy the {\em Feller property}
(i.e. corresponding Markov operators leave the space of bounded
continuous functions invariant) or some its generalizations
and the compactness of the phase psace.
The situation is rather different if one wants to study
the existence of an invariant measure supported by a
non compact subset of the phase space as in the case
considered in Section~\ref{s:inv-mes-density}.

We follow the convention: the upper index is always reserved for
time, while the lower index is used for space variables.
The only exception is the upper index of the shift map
$\sigma^\ell$ which means the spatial translation by the
vector $\ell$ and, in fact, can be interpreted as a time
variable having in mind that we apply $|\ell|$ times the
translation by the length one in the direction $\ell$.

\section{Functional-analytic approach}\label{s:functional}

For a Markov chain $\xi^t$ on a metric space $(X,\dist)$ denote by %
$$\P_x(\xi^t\in{A}):=P^t(x,A), \quad
  \E_x(\phi(\xi^t)):=P^t\phi(x) $$ %
the corresponding probability and mathematical expectation for
the trajectory started from the point $x\in X$. Recall that a
measure $\mu$ is invariant (stationary) for the Markov chain
$\xi^t$ if
$$ \mu(\E_x(\phi(\xi^t))) = \mu(\phi) $$
for any continuous function $\phi:X\to\IR$, where
$\mu(\phi):=\int_X\phi~d\mu$.

We start with a well known folklore conditional result about the
uniqueness of an invariant measure of a Markov process assuming
that this measure does exist.

\begin{theorem}\label{t:uniqueness} Assume that there exists a
probability invariant measure $\mu$ for a Markov chain $\xi^t$
and let %
\beq{e:intersection}{\P_{x,y}(\tau^0<\infty)=1 ~\quad\forall x,y\in X.}%
Then $\mu$ is the only probability invariant measure of this
Markov chain and $P^t_*\nu$ converges in the total variation
metric to $\mu$ as $t\to\infty$ for any probability measure
$\nu$.
\end{theorem}

\proof Consider a {\em coupled} system $(\xi^t,\acute\xi^t)$ whose
trajectories are assumed to coincide after their first
intersection. Integrating the inequality (\ref{e:coupling-ineq})
with respect the measures $\mu$ and $\nu$ we have %
$$ |P_*^t\nu(A) - \mu(A)|
\le \int\!\!\!\!\int\P_{x,y}(\tau^0>t)~d\nu(x)~d\mu(y) $$ %
which vanishes as $t\to\infty$ due to the assumption
(\ref{e:intersection}).\qed

In a more general setting different realizations of a Markov chain
may not intersect but only converge to each other in time. It
turns out that the previous result can be generalized to this
setting as well (at least partially).

\begin{theorem}\label{t:uniqueness-weak} Assume that there exists a
probability invariant measure $\mu$ for a Markov chain $\xi^t$ and let %
\beq{e:intersection-weak}{\P_{x,y}(\tau^\ep<\infty)=1
                          \quad\forall x,y\in X, ~\ep>0.}%
Then $\mu$ is the only probability invariant measure of this
Markov chain and $P^t_*\nu\toas{t\to\infty}\mu$ weakly
for any probability measure $\nu$.
\end{theorem}

\proof Recall that a function $\phi:X\to\IR$ is called Lipschitz
continuous if there exists a finite constant ${\rm Lip}(\phi)$
such that %
\beq{e:Lip}{|\phi(x) - \phi(y)|\le{\rm Lip}(\phi)\dist(x,y) }%
for any $x,y\in X$.

Define $|\phi|_L:=|\phi|_\infty+{\rm Lip}(\phi)$ and consider only
functions $\phi$ with bounded $|\cdot|_L$-norm.
We have: %
$$|\E_x\phi(\xi^t) -\E_y\phi(\acute\xi^t)| %
= |\E_{x,y}\phi(\xi^t) -\E_{x,y}\phi(\acute\xi^t)| %
\le \E_{x,y}|\phi(\xi^t) - \phi(\acute\xi^t)| \CR%
\le 2|\phi|_\infty \P_{x,y}(\tau^\ep\ge t)
  + \ep{\rm Lip}(\phi) \P_{x,y}(\tau^\ep<t)%
\toas{t\to\infty}\ep{\rm Lip}(\phi).$$ %
Now since the left hand side does not depend on
$\ep>0$ we deduce that
$|\E_x\phi(\xi^t) -\E_y\phi(\acute\xi^t)|\toas{t\to\infty}0$.

On the other hand, for any measure $\nu$ %
\beq{e:eqq}{\mu(\phi)-P_*^t\nu(\phi)
 = (\mu\times\nu)(\E_x\phi(\xi^t) -\E_y\phi(\acute\xi^t)) .}%

Observe that on a locally compact space the Lipschitz continuous
functions are dense in the set of all continuous functions being
constant at infinity (by Stone-Weierstrass Theorem). So to
distinguish between any two probability measures on such space
it is enough to consider only integrals of Lipschitz continuous
functions. Hence the invariant measure is unique and 
$P^t_*\nu\toas{t\to\infty}\mu$ weakly. \qed

A close look at the proof shows that in the relation (\ref{e:eqq}) 
one needs to control the integrands only up to the sets of zero 
measure. Therefore the uniformity of the 
assumption~(\ref{e:intersection-weak}) on the choice of pairs of 
initial conditions $(x,y)$ may be significantly relaxed:

\begin{corollary}\label{c:intersection-weak} 
Let $\mu,\nu$ be probability invariant measures for a Markov chain 
$\xi^t$ and let %
\beq{e:intersection-weak-weak}{\P_{x,y}(\tau^\ep<\infty)=1
                          \quad\forall \ep>0}%
for $\mu\times\nu$ almost all pairs of points $x,y\in X$. 
Then $\mu=\nu$.
\end{corollary}

\section{Probability techniques}\label{s:probability}%
In this section we shall study the topological quasi successful
couplings.

\begin{theorem}\label{t:uniq-prob} Let $\xi^t$ be a Markov chain
with an invariant measure $\mu$ and let $\acute\xi^t$ be another
version of the same Markov chain. Assume that there exists a
topologically quasi successful coupling for the Markov chains
$\xi^t,\acute\xi^t$, i.e. the relation~(\ref{e:t-q-success})
holds. Then $\mu$ is the only probability invariant measure of
this Markov chain and $P^t_*\nu\toas{t\to\infty}\mu$ weakly
for any probability measure $\nu$.
\end{theorem}

\proof For an open set $A\in\cB$ and a pair of points
$x,y\in X$ we introduce the notation %
\beq{e:phi}{\Phi^{t}(x,y,A):=\
            |\P_{x,y}(\xi^t\in A)-\P_{x,y}(\acute\xi^t\in A)|.}%
Then %
$$ \P_{x,y}(\xi^t\in A)
 = \P_{x,y}(\xi^t\in{A},\acute\xi^t\in{A})
 + \P_{x,y}(\xi^t\in{A},\acute\xi^t\notin{A})$$
and
$$ \P_{x,y}(\acute\xi^t\in A)
 = \P_{x,y}(\xi^t\in{A},\acute\xi^t\in{A})
 + \P_{x,y}(\xi^t\notin{A},\acute\xi^t\in{A}).$$
Hence, %
\bea{ \Phi^{t}(x,y,A)
 \a= \left|\P_{x,y}(\xi^t\in{A},\acute\xi^t\notin{A})
       - \P_{x,y}(\xi^t\notin{A},\acute\xi^t\in{A})\right| \\
 \a\le \P_{x,y}(1_A(\xi^t)\ne1_A(\acute\xi^t))
   \toas{t\to\infty}0 , } %
by the assumption~(\ref{e:t-q-success}).

On the other hand, since %
$$\P_{x,y}(\xi^{t}\in A)=\P_{x}(\xi^{t}\in A)=P^{t}(x,A)$$
and
$$\P_{x,y}(\acute\xi^{t}\in A) %
 =\P_{y}(\acute\xi^{t}\in A)=P^{t}(y,A)$$
from the relation~(\ref{e:phi}) we get %
$$ -\Phi^{t}(x,y,A) \le P^{t}(x,A) - P^{t}(y,A) \le \Phi^{t}(x,y,A) .$$
Integrating the last inequality with respect to the probability
invariant measure $d\mu(x)$ and using that %
$$ \int P^{t}(x,A)~d\mu(x) = P_{*}^{t}\mu(A)=\mu(A) $$
we get
$$ -\int\Phi^{t}(x,y,A)~d\mu(x) \le \mu(A) - P^{t}(y,A)
   \le \int\Phi^{t}(x,y,A)~d\mu(x) .$$
Thus %
$$ |\mu(A) - P^{t}(y,A)| \le \int\Phi^{t}(x,y,A)~d\mu(x) .$$
Therefore $P^{t}(y,A)\toas{t\to\infty}\mu(A)$. Any bounded continuous 
function $\phi$ can be approximated uniformly by a finite combination 
of indicator functions. 
Therefore we have %
\bea{ P_{*}^{t}\nu(\phi) \a= \nu(P^{t}\phi)
                      = \nu(\int \phi(z)P^{t}(\cdot,dz)) \\
                      \a\toas{t\to\infty} \nu(\mu(\phi))
                      =\mu(\phi)\nu(1)=\mu(\phi)  }%
which implies the weak convergence of $P_{*}^{t}\nu$. \qed

A close look at the above proof shows that if the
relation~(\ref{e:t-q-success}) holds for a given set
$A\in\cB$ and any $x,y\in X$ one has
$$ P_{*}^{t}\nu(A)\toas{t\to\infty}\mu(A) $$
for any probability measure $\nu$. Thus at least on the
set $A$ the limit measure coincides with $\mu$.

The following simple generalization of this result is useful for
the analysis of interacting particle systems. Let $S:X\to X$ be
a measurable bijection such that $P^{t}(x,A)=P^{t}(Sx,SA)$ for any
$x\in X,~A\in\cB,~t\ge0$.

\begin{corollary}\label{c:prob} Assume that all conditions of
Theorem~\ref{t:uniq-prob} hold except that instead of the
topologically quasi successful coupling condition one assumes that
there exists a coupling such that %
\beq{e:t-q-success-shift}{
     \P_{x,y}(1_A(\xi^t)\ne1_A(S\acute\xi^t))
     \toas{t\to\infty}0  \quad\forall A\in\cB, ~x,y\in X,}%
Then the claim of Theorem~\ref{t:uniq-prob} remains valid.
\end{corollary}

The proof of this result follows exactly the same lines as above
except that one uses $S\acute\xi^{t}$ instead of $\acute\xi^{t}$
and later $Sy$ instead of $y$ since
$$ \P_{y}(S\acute\xi^{t}\in A) = \P_{y}(\acute\xi^{t}\in S^{-1}A)
 = P^{t}(y,S^{-1}A) = P^{t}(Sy,A) .$$

Similarly to Corollary~\ref{c:intersection-weak} the result of 
Theorem~\ref{t:uniq-prob} can be generalized for the case when 
the relation~(\ref{e:t-q-success}) holds only for almost all 
pairs of initial points $x,y\in X$.

\section{Markovian couplings}\label{s:markov}

A coupling is called {\em Markovian} if the pairs process
$(\xi^{t},\acute\xi^{t})$ is a Markov chain on the direct product
space $X\times X$. Under this assumption the proofs of our
results about the uniqueness of the invariant measure might be
significantly simplified by means of the following idea borrowed
from \cite{Lig,EG}. Observe that if the coupling is Markovian and
the Feller Markov chain $\xi^t$ admits two probabilistic invariant
measures $\mu\ne\nu$ then the Markov pairs process
$(\xi^{t},\acute\xi^{t})$ admits an invariant measure having
marginals equal to $\mu$ and $\nu$ correspondingly. Assume now
that any version of the quasi  successful coupling takes place.
Then under the coupling the processes $\xi^t$ and $\acute\xi^t$
eventually behave the same, which contradicts to the assumption
that $\mu\ne\nu$.

Let $(\xi^{t},\acute\xi^{t})$ be a Markovian coupling of two
copies of the Markov chain $\xi^{t}$. Denote by $\cM^{t}$ the
distribution of the pair $\xi^{t},\acute\xi^{t}$ at time $t\ge0$
and by $\pi_{*}\cM^{t}$ and $\acute\pi_{*}\cM^{t}$ the projections
of the measure $\cM^{t}$ to the components of the coupled system.
Then from a measure-theoretic point of view both our quasi
successful coupling assumptions mean the the projections
$\pi_{*}\cM^{t}$ and $\acute\pi_{*}\cM^{t}$ are becoming close to
each other in a certain sense. On the other hand, the convergence
of the projections $\pi_{*}\cM^{t}$ and $\acute\pi_{*}\cM^{t}$
immediately implies that under the action of the original Markov
operator $P_{*}^{t}$ any two probability measures converge to each
other.

The Markovian assumption is very natural and, in fact, is
satisfied in all main coupling constructions. Using the idea
proposed in \cite{Ros} under this assumption one can improve
significantly the result of Theorem~\ref{t:uniqueness}.

\begin{theorem}\label{t:uniqueness-f} Let $(\xi^{t},\acute\xi^{t})$
be a Markovian coupling of two copies of a Markov chain on the
probability space $(\Omega,\P_{x,y})$ and let $\tau$ be a random
variable such that $\xi^{\tau}=\acute\xi^{\tau}$ and
$\P_{x,y}(\tau<\infty)=1$. Then if there exists a probability invariant
measure $\mu$ of the Markov chain $\xi^{t}$ then $P^t_*\nu$
converges in the total variation metric to $\mu$ as $t\to\infty$
for any probability measure $\nu$.
\end{theorem}
\proof Define a new process: %
$$ \t\xi^{t} := \function{\acute\xi^{t} &\mbox{if } t\le\tau\\
                          \xi^{t}       &\mbox{otherwise} .} $$
In \cite{Ros} it has been shown that if the coupling is Markovian
then the process $\t\xi^{t}$ is Markov with the same transition
probabilities $P^{t}(\cdot,\cdot)$. Thus the coupling
$(\xi^{t},\t\xi^{t})$ is successful according to our definition
and the result follows from Theorem~\ref{t:uniqueness}. \qed

One is tempted to use a similar argument to simplify the
assumptions about the quasi successful couplings. Unfortunately
this does not work even in the deterministic setting. Consider a
binary map $Tx:=2x~({\rm mod} 1)$ defined on the unit interval
$X:=[0,1]$. Then the direct product of these maps is a coupled
process with the Markovian coupling. Moreover, for Lebesgue almost
any initial point $x\in X$ and any $y\in X$ the trajectories 
$\xi^{t}:=T^{t}x$
and $\acute\xi^{t}:=T^{t}y$ are becoming arbitrary close
infinitely many times. On the other hand, the dynamical system
$(T,X)$ has infinitely many ergodic probability invariant
measures, in particular, the Lebesgue measure and the
delta-measure at zero, which contradicts to 
Corollary~\ref{c:intersection-weak}.

\section{Applications to particle systems}\label{s:particles}
By discrete space infinite particle systems one means the
translational invariant dynamics of infinite configurations of
particles on an integer lattice. We refer the reader to \cite{Lig}
for a general discussion of systems of this sort.

To apply the coupling technique to discrete space particle
dynamics one needs to construct a coupling between different
realizations of the dynamics. In distinction to the finite
dimensional systems there is no much hope to achieve the usual
successful coupling since this would mean that an infinite number
of pairs of particles in the components of the coupled system need
to share the same positions simultaneously. It seems more probable
to get a version of a quasi successful coupling but even this
setting turns out to be restrictive. The only known way to
overcome this difficulty is to use a kind of a {\em consecutive
coupling} based on pairing. By a {\em pairing} of two particles
belonging to two different configurations we shall mean that up
to a certain moment of time, which we call the {\em pairing
time}, these particles behave independently of each other and
after that time the rules applied to these two particles in both
processes are assumed to be the same if this does not contradict
to the definition of the process. If the latter event happen the
pair needs to be broken.

The most known construction, which we shall call an {\em equal}
pairing is based on a very simple idea that two particles
belonging to different configurations are becoming paired when
they are located simultaneously at the same lattice site. This
construction has been successfully applied e.g. in \cite{Lig} for
a continuous time totally asymmetric exclusion process (TASEP).
However even for a discrete time version of the TASEP (not
speaking about more general exclusion type processes) the equal
coupling does not work well since in the discrete time setting
several particles may move simultaneously (which cannot happen in
the continuous time case). Another also consecutive coupling
making use of a non equal pairing has been proposed to solve this
problem for the discrete time TASEP by L.~Gray \cite{Gr}, see
also our forthcoming paper for the analysis of more
general consecutive coupling constructions.
Roughly speaking instead of pairing of particles sharing the same
positions L.~Gray proposed to pair particles located close enough.
The main difference to the equal pairing is that one proves that
eventually with time only pairs of particles with the
same distances between their members survive almost surely.

Even in the case of the equal pairing one cannot apply directly
our results about quasi successful couplings since their
conditions are still too strong. Typically one has to deal with
one of the following scenarios:

\begin{enumerate}
\item the upper density of the set of discrepancies between
particle configurations (possibly calculated up to a finite
spatial shift) vanishes with time;

\item for any finite segment $I\in\IZ^{d}$ the probability that
particle configurations do not coincide on this segment (possibly
calculated after a finite spatial shift) vanishes with time.
\end{enumerate}

Without spatial shifts these scenarios can be considered as direct
generalizations of the quasi successful coupling and the weak
quasi successful coupling respectively. Much weaker versions with
spatial shifts are necessary to take into account situations when
non equal pairs of particles survive in a coupled system. Note
that despite both scenarios look very similar none of them is a
consequence of another one. Indeed the presence of a finite number
of fixed discrepancies does not change their density while it
changes completely the situation described by the second scenario.

Specific models of particle systems and specific pairing for them
leading to above conditions will be studied in a separate paper
and here we restrict ourselves to ergodic consequences
of these scenarios.

\subsection{Functional-analytic setting}
Let $X:=\cA^{\IZ^{d}},~d\ge1$ be a space of be-infinite sequences
with elements from a finite alphabet $\cA$ endowed with the
standard metric $\dist(\cdot,\cdot)$ based on cylinders, namely %
$$ \dist(x,y) := 2^{-\kappa(x,y)} ,$$
where %
$x,y\in X$, $\kappa(x,y)$ is defined as the largest positive
integer $n$ such that $x_\ell=y_\ell$ for all $\ell\in\IZ^{d}$
such that $|\ell|\le{n}$, and $|\ell|:=\max_{i}\ell_{i}$.

It is straightforward to check that the metric space $(X,\dist)$
is compact. Recall that two probability measures $\mu,\nu$ are
equal if $\mu(\phi)=\nu(\phi)$ for any continuous function
$\phi:X\to\IR$ and different otherwise. In this sense measures
can be distinguished by means of continuous functions.

To be consistent with the previous results we fix also the
standard sigma-algebra of Borel sets $\cB$ on $X$ based on
cylindric sets and denote by $\cM$ the set of all translational 
invariant probability measures $m$ on $X$. The triple $(X,\cB,\dist)$
forms a compact measurable metric space on which we consider our
Markov chains.

Since Lipschitz continuous functions are dense in the set of all
continuous functions it is enough to use only Lipschitz continuous
test functions to distinguish between any two distinct measures.

We introduce a semimetric $\tilde\dist(x,y)$ equal to the upper
density of the set of discrepancies between the configurations
$x,y\in X$ as follows: %
$$ \tilde\dist(x,y) := \limsup_{n\to\infty}
         \frac{\#\{k\in\IZ^{d}: ~~x_k\ne y_k,~|k|\le{n}\}}{(2n+1)^{d}} $$
and say that a function $\phi:X\to\IR$ is $\tilde\dist$-Lipschitz
continuous if there exists a finite constant %
$\widetilde{\rm Lip}(\phi)$ such that the inequality %
\beq{e:Lip-semi}{ |\phi(x) - \phi(y)| \le %
  \widetilde{\rm Lip}(\phi)\tilde\dist(x,y) }%
holds almost everywhere (a.e.) with respect to any translationally 
invariant measure.

Using this notation the first scenario mentioned above can be
reduced to the question up to which extent measures can be
distinguished by means of $\tilde\dist$-Lipschitz test functions.

Clearly there are measures %
which cannot be distinguished this way. Therefore we restrict this
question to translational invariant measures. Moreover, we shall
restrict ourselves to some subset $X_r$ of the set of 
configurations invariant with respect to dynamics and to spatial 
shifts (e.g. the set of configurations of a given particle 
density). Denote by $\sigma^\ell\xi^t$ the spatial shift of the 
configuration$\xi^t$ by the integer vector $\ell\in\IZ^{d}$, i.e. 
$(\sigma^\ell\xi^t)_i:=\xi^t_{i+\ell}$.

\begin{theorem}\label{t:uniq-anal'} Let there exist a coupling of
two identical particle systems $\xi^t,\acute\xi^t$ such that for
any two initial configurations $\xi^0,\acute\xi^0\in X_r$ we have %
\beq{e:equal-anal} {
  \tilde\dist(\xi^t,\sigma^\ell\acute\xi^t)\toas{t\to\infty}0 }%
for some $\ell\in\IZ^d,~|\ell|<L$ with probability one.
Then the cardinality of the set of translational (both in space
and time) invariant probability measures of the Markov chain
$\xi^t$ supported by $X_r$ cannot exceed one.
\end{theorem}
\proof We start with the stronger version of the
assumption~(\ref{e:equal-anal}) without the spatial shifts namely
we assume that: %
\beq{e:equal-anal1} {\tilde\dist(\xi^t,\acute\xi^t)\toas{t\to\infty}0 }%
with probability one.

Let $\mu,\nu$ be two different translational, both in
space and time, invariant probability measures of the particle
system. Then there exists a continuous function $\phi$ which can
distinguish between them, i.e. $\mu(\phi)\ne\nu(\phi)$. On the
other hand, any continuous function can be approximated by
indicators of finite cylinders. Thus there exists a finite
cylinder $\cyl{I}$ with the base $I\subset\IZ^d$ for which
$\mu(\cyl{I})\ne\nu(\cyl{I})$. Recall that a cylinder with the
base $I$ is the set of configurations $x\in X$ with the given
subset of coordinates $x_i$ for $i\in I$, i.e.
$\cyl{I}:=\{x\in X: ~~x_i=a_i,~i\in I\}$.

Consider a sequence of functions %
$$ \psi_n:=\frac1{(2n+1)^{d}}\sum_{|\ell|\le n}
                 1_{\cyl{I}}\circ\sigma^\ell ,$$
where $1_A(\cdot)$ is the indicator function of the set $A\in X$.

Observe that the function $\psi_0:=1_{\cyl{I}}$ distinguishes
between the measures $\mu,\nu$. Therefore due to the translation
invariance of these measures each of the functions $\psi_n(\cdot)$
again can distinguish between the measures $\mu,\nu$. Moreover,
$$ \mu(\psi_n)
 = \frac1{(2n+1)^{d}}\sum_{|\ell|\le n} 1_{\cyl{I}}\circ\sigma^\ell
 = \frac1{(2n+1)^{d}}\sum_{|\ell|\le n} \mu(1_{\cyl{I}})
 = \mu(\psi_0) $$
for any $n\in\IZ_+$ and thus
$$ |\mu(\psi_n) - \nu(\psi_n)|
 = |\mu(\psi_0) - \nu(\psi_0)|
 = |\mu(1_{\cyl{I}}) - \nu(1_{\cyl{I}})| >0.$$
The limit function $\psi(x):=\lim_{n\to\infty}\psi_n(x)$ exists a.e.
(with respect to any translational invariant measure) by the
Birkhoff Ergodic Theorem. Due to the estimates above the limit
function again distinguishes between the measures $\mu,\nu$.

Let us show that this limit function $\psi$ is $\tilde\dist$-Lipschitz
continuous. For $x\in X$ and $n\in\IZ_+$ denote
$$ J_n(x):=\{\ell\in\IZ^d:~~ |\ell|\le n,
                          ~~ \sigma^\ell(x)\in\cyl{I}\}.$$
Then $\psi_n(x) = (2n+1)^{-d}|J_n(x)|$. A single site may
contribute to the value $|J_n(x)|$ at most a certain finite
constant $C$ which depends only on the cylinder $\cyl{I}$
but neither on $n\in\IZ_+$ nor on the configuration $x\in X$.

For a pair of configurations $x,y\in X$ and $n\in\IZ_+$ denote
$$ D_n(x,y):=\{k\in\IZ^{d}: ~~x_k\ne y_k,~|k|\le{n}\} .$$
Then
$$ |\psi_n(x) - \psi_n(y)| \le C (2n+1)^{-d}~|D_{n+|I|}(x,y)| .$$
On the other hand,
$$ \limsup_{n\to\infty} (2n+1)^{-d}~|D_{n}(x,y)|
 = \tilde\dist_n(x,y) .$$
Therefore passing to the limit as $n\to\infty$ we get that
the limit function $\psi(\cdot)$ is $\tilde\dist$-Lipschitz
continuous with the constant $\widetilde{\rm Lip}(\psi):=C$.

Applying now the argument used in the first part\footnote{
   Since we deal here only with Lipschitz continuous functions
   we do not need to apply the last part of the proof of
   Theorem~\ref{t:uniqueness-weak} which deals with the
   approximations of continuous functions by the Lipschitz
   continuous ones.} %
of the proof of Theorem~\ref{t:uniqueness-weak} we get the result.

Using similar arguments one considers also the case when the
components of the coupled process converge to each other in terms
of the semimetric $\tilde\dist$ only after a finite spatial shift.
\qed

\subsection{Probabilistic setting}

Assume now that there exist a coupling of two identical particle
systems such that for any given finite subset $I\in\IZ$ the
restriction of the components of the coupled system to $I$
coincide at moment $t\ge0$ with probability going to one with
time, i.e. %
$$ \P_{x,y}(\xi^t_{|I}=\acute\xi^t_{|I})\toas{t\to\infty}1 $$
for any two initial configurations $\xi^0,\acute\xi^0$ from a
certain dynamically invariant subset $X_r\subseteq X$. Similar
to the deterministic setting we say that a set $Y\subseteq X$ is
{\em dynamically invariant} if $\xi^0\in Y$ implies that
$\xi^t\in Y$ a.s. for any $t\ge0$.

Thus the probability that the configurations $\xi^t,\acute\xi^t$
coincide on any finite cylinder goes to one with time. Hence the
relation~(\ref{e:phi}) is satisfied for any cylinder set. Using
now that linear combinations of indicator functions of the cylinder sets 
are dense in the space of bounded continuous functions on $X$ we 
may apply the claim of Theorem~\ref{t:uniq-prob}.

In the case when the restriction of the configurations to a finite
segment $I$ coincide only after a finite spatial shift
$\sigma^{\ell}$ follows from the same argument if we consider only
(spatially) translational invariant distributions and use
Corollary~\ref{c:prob} with $S:=\sigma^{\ell}$ instead of
Theorem~\ref{t:uniq-prob}.

Therefore we get the following result. Denote by $\cM_{r}$ the set
of translational invariant probability measures supported by
$X_{r}$.

\begin{theorem}\label{t:uniq-prob'} Let $\xi^t$ be a Markov chain
with an invariant measure $\mu\in\cM_{r}$ and let $\acute\xi^t$
be another version of the same Markov chain. Assume that there
exists a coupling under which
$$ \P_{x,y}(\xi^t_{|I}
  =\sigma^\ell\acute\xi^t_{|I})\toas{t\to\infty}1 $$ %
for any pair of configurations $\xi^0=x,~\acute\xi^0=y\in X_{r}$,
any finite segment $I\subset\IZ^{d}$, and some
$\ell\in\IZ^d, ~|\ell|\le L<\infty$. Then $P^t_*\nu$ converges
weakly to $\mu$ as $t\to\infty$ for any $\nu\in\cM_{r}$.
\end{theorem}

\begin{corollary} Under all the assumption of the previous Theorem 
except for $\mu\in\cM_{r}$ the sequence $P^t_*\nu$ converges
weakly to some invariant measure which does not depend on $\nu$. 
\end{corollary}

\section{Invariant measures with a given particle density}
\label{s:inv-mes-density}

As we already mentioned normally the proof of the existence
of invariant measures can be reduced to the checking of some
very general topological and metrical assumptions. The
situation changes drustically if one wants to study invariant
measures having some additional properties, e.g. supported
by some dynamically invariant subsets. The problem here is that
even if the phase space of the system is compact the dynamically
invariant subset may not satisfy this property.

To be more specific consider a class of locally interacting
particle systems on the lattice $\IZ^d$ represented by Markov
chains with the phase space 
$X:=\cA^{\IZ^d}, ~ \cA:=\{0,1,\dots,|\cA|-1\}, ~|\cA|<\infty$.
For a configuration $x\in X$ and a lattice site $i\in\IZ^d$ we
associate the value $x_i\in\cA$ to the number of particles located
at $i$. The dynamics is defined as follows. At each site of the
lattice there is an alarm-clock and at time $t>0$ we consider
only those lattice sites at which the alarm rings. For each particle
located at such site $i\in\IZ^d$ we calculate its {\em velocity}
$v_i$, such that $|v_i|\le V<\infty$, using a (random or
deterministic) procedure which does not depend neither on time nor
on the configuration $x$. Here $V$ plays the role of the largest
allowed velocity and its boundedness defines the locality of
interactions. Then one checks a certain {\em admissibility}
condition related to the possibility to move a particle from the
site $i$ to the site $i+v_i$. We assume that the admissibility
condition is again local and depends only on the present positions
of the particles in the $2V$-lattice neighborhood of the site $i$.
Only if this condition is satisfied the particle is moved to
a new position. Then for all sites to where the particles were
moved we restart the alarm-clocks (again using a certain random
or deterministic procedure).

Depending on the way how one restarts the alarm-clocks both
continuous and discrete time particle systems can be considered.
In what follows we restrict ourselves to a (more interesting
from our point of view and much less studied) discrete time
case, assuming that the alarm-clocks start with the same setting
and after each restart we add one to the time. Therefore all
particles are trying to move simultaneously.

There is an important property that holds for all systems we
consider here: particle number conservation. For a configuration
$x\in X:=\cA^{\IZ^d}$ and a finite subset $I\subset\IZ^d$
denote by $\den(x,I)$ the number of particles from the
configuration $x$ located in $I$ divided by the total number
of sites in $I$, which we denote by $|I|$. Clearly
$0\le\den(x,I)\le|\cA|-1$. Choosing a sequence of lattice cubes
$I_n:=\{\ell\in\IZ^d:~|\ell|\le n\}$ we consider the limit
$\lim_{n\to\infty}\den(x,I_n)$. If this limit exists we call
it the {\em particle density} of the configuration $x\in X$
and denote by $\den(x)$.

According to the above informal description the particle
systems under consideration are Markov chains on $X$. We assume
the they satisfy the Feller condition, which together with
the compactness of the phase space implies immediately
the existence of an invariant measure. Moreover, due to the
particle number conservation it seems that for each
$r\in[0,|\cA|-1]$ there should be an invariant measure
$\mu_r$ having a support on the set of configurations of
density $r$.

The following simple deterministic example shows that this
is absolutely not the case. For $d=1$ define a map
$T:X\to X$ as follows:
$$  (Tx)_i:=\function{x_i     &\mbox{if } i=0 \\
                      0       &\mbox{if } |i|=1 \\
                      x_{i-1} &\mbox{if } i>1 \\
                      x_{i+1} &\mbox{otherwise}.} $$
In words, we preserve the central site, shift all
other sites to the left and to the right, and set the
sites $\pm1$ to zero.

Clearly all above assumptions are satisfied in this example and
the particle density is conserved. On the other hand, for any
initial configuration $x\in X$ the sequence of configurations
$T^tx$ converges in the standard metric to the configuration of
all zeros, having, of course, zero particle density. Thus this
system has only one invariant measure supported by the
configuration of all zeros.

One might think that the result above is a consequence of the
absence of the conservation of vacancies. Let us describe a simple
probabilistic generalization of this example conserving both
particles and vacancies and having exactly two invariant measures
supported by the configuration of all zeros or all ones. Choose
$p\in(0,1)$. For each particle in a configuration $x\in X$ located
at a nonnegative site and having a vacancy immediately to the
right we exchange the positions of these particle/vacancy with the
probability $p$. If the particle is located at a negative site we
do the same but for the vacancy immediately to the left of it.
This rule defines a discrete time Markov chain having the desired
property.

\bigskip

In what follows we discuss a coupling argument\footnote{The
   idea of this construction has been proposed by L.~Gray
   \cite{Gr} for the case of the one-dimensional totally
   asymmetric exclusion process. It is worth note that
   the situation we consider differs significantly from the
   one studied by L.~Gray in that particles may have
   velocities greater than one and thus long range interactions
   should be taken into account.} %
which gives sufficient conditions for the existence of invariant
measures supported by configurations of a given density.

Let $\xi^t$ be a Markov chain on $X:=\cA^{\IZ^d}$ with
translational invariant transition probabilities,
i.e. $P(x,A)\equiv P(\sigma^\ell{x},\sigma^\ell{A})$
for any $x\in X, ~A\in\cB, ~\ell\in\IZ^d$, which describes
the dynamics of a locally interacting particle system.
We shall say that a {\em particle density} $\den(\mu)$ of
a translational invariant probability measure $\mu$ is the
mathematical expectation (with respect to this measure) of
the number of particles located at the origin,
i.e. $\den(\mu):=\mu(x_0)$.

The following simple result explains the reason why we use
the same notation for the the particle densities of
configurations and measures.

\begin{lemma}\label{l:transl-inv-ex} For each $r\in[0,|\cA|-1]$
there exists a translational both in space and time invariant
measure $\mu_r$ such that $\den(\mu_r)=r$.
\end{lemma}
\proof Choose any translational invariant measure $\mu$ and
consider a sequence of measures
$$ \mu^n:=\frac1n \sum_{t=0}^{n-1} P_*^t\mu .$$
Due to the compactness of $X$ there exists a subsequence
$n_k\to\infty$ such that $\mu^{n_k}$ converges weakly
to a limit measure $\nu$. Now using the Feller property
we see that the measure $\nu$ is invariant. The translation
invariance of $\nu$ follows from the construction and
since the particle density of each measure $\mu^n$ is the
same as for the initial measure $\mu$, we get
$\den(\nu)=\den(\mu)$. Consider a distribution
$\pi:=\{\pi_i\}_{i\in\cA}$ such that $\sum_i i\pi_i=r$.
To finish the proof it remains to notice that the Bernoulli
measure with the alphabet $\cA$ and the distribution of
symbols $\pi$ has the particle density which is exactly
equal to $r$. \qed

\begin{corollary}\label{c:bernoulli} For each $r\in[0,|\cA|-1]$ 
the measure $\mu_r$ with properties described above can be constructed as 
a limit point of the sequence of ergodic averages of the measures 
$P_*^t\mu_{r,B}$ where $\mu_{r,B}$ is a Bernoulli measure with the 
particle density $r$. 
\end{corollary}

Unfortunately a translational invariant measure with the
particle density $r$ is not necessarily supported by
configurations of density $r$.

We shall say that a measure $\mu$ is {\em spatially
ergodic} if the dynamical system $(\sigma, X,\mu)$ is
ergodic. We refer the reader to \cite{Sinai} for general
definitions related to ergodic theory of dynamical systems.

Clearly each spatially ergodic measure is translational
invariant and by the Birkhoff theorem it is supported
by the configurations of the same particle density.

\begin{lemma}\label{l:desity-preservation} Let the initial
distribution $\mu$ of the process $\xi^t$ be spatially
ergodic. Then almost surely for any $t\ge0$ the density
$\den(\xi^t)$ for the configuration $\xi^t$ is well defined,
does not depend on $t$, and is equal to the particle
density $\den(\mu)$ of the initial distribution.
\end{lemma}
\proof For $t=0$ the existence of the limit in the definition
of the particle density for the configuration $\xi^0$ and
its coinscidence with $\den(\mu)$ follows from the Birkhoff
ergodic theorem due to the spatial ergodicity of the initial
distribution $\mu$.

Choose $n\in\IZ_+$ large enough and consider the lattice
cube $I_n$. We can estimate the difference between
the number of particles located in $I_n$ at moments $t$
and $t+1$ as follows. Only particles belonging to the
set $I_n\setminus I_{n-V}$ are able to leave the
the cube $I_n$ during one time step. Similarly only
particles located at the set %
$I_{n+V}\setminus I_{n}$ are able to enter the cube
$I_{n}$ during this time. Therefore using that
$0\le\den(\xi^t,I)\le|\cA|-1$ we get %
\bea{ \a|\den(\xi^{t+1},I_n) - \den(\xi^t,I_n)| \\
 \a\le \frac{|I_n\setminus I_{n-V}|}{|I_n|}
       ~\den(\xi^t,I_n\setminus I_{n-V})
     ~~+~~ \frac{|I_{n+V}\setminus I_n|}{|I_n|}
       ~\den(\xi^t,I_{n+V}\setminus I_n)  \\
 \a\le 2|\cA|\max\left(
        \frac{(2n+1)^d - (2n+1-2V)^d}{(2n+1)^d}, ~~
        \frac{(2n+1+2V)^d - (2n+1)^d}{(2n+1)^d}
                 \right) \\
 \a= 2|\cA|\max\left(1-\left(1-\frac{2V}{2n+1}\right)^d, ~~
                \left(1+\frac{2V}{2n+1}\right)^d-1 \right)
     \toas{n\to\infty}0 .}%
This proves both the existence and the conservation of the
particle density for the configurations $\xi^t$. \qed

\begin{lemma}\label{l:transl-inv-erg}
Let there exist a coupling of two copies
$\xi^t,\acute\xi^t$ of our Markov chain such that for all
initial configurations with the same
particle density the conditions of either
Theorem~\ref{t:uniq-anal'} or Theorem~\ref{t:uniq-prob'}
hold. Then for each spatially ergodic initial distribution
$\mu$ all limit points of the sequence of measures
$\{P_*^n\mu\}_{n\in\IZ_+}$ coincide and the limit measure
has the same particle density.
\end{lemma}
\proof Fix a spatially ergodic distribution $\mu$ on $X$.
Due to the compactness of $X$ there exists a sequence of
integers $n_k\toas{k\to\infty}\infty$ such that the
sequence of measures $P_*^{n_k}\mu$ converges to a
certain limit point, call it $\t\mu$. Passing to the
limit as $n_k\to\infty$ and using that $\den(P_*^{n_k}\mu)$
does not depend on $n_k\in\IZ_+$ we get that the limit
distribution $\t\mu$ also has the same particle density. 
The spatial ergodicity of $\mu$ implies the spatial ergodicity
of $P_*^n\mu$ for each $n\in\IZ_+$.

Coupling\footnote{Observe that to make our coupling between
   two copies of the process with certain initial distributions
   we need those distributions to be spatially ergodic and to
   apply any kind of a quasi successful property they need to
   be of the same particle density.} %
the process with the initial distribution $\mu$ and the one
with $P_*\mu$ we see that according to our assumptions
these two processes eventually behave the same. Using the 
same arguments as in the proofs of Theorems~\ref{t:uniqueness-weak} 
and \ref{t:uniq-prob} we get $P_*^{n_k+1}\mu\toas{n_k\to\infty}\t\mu$.
Thus using the semigroup property
$P_*^{n+1}=P_*^nP_*$ we have:
$$ \t\mu = \lim_{k\to\infty} P_*^{n_k+1}\mu
         = P_* \left( \lim_{k\to\infty}P_*^{n_k}\mu \right)
         = P_*\t\mu ,$$
which proves the invariance of the measure $\t\mu$, and
in turn implies that all limit points of the sequence of
measures $P_*^n\mu$ coincide and do not depend on $\mu$. \qed

Finally let us discuss briefly a prototype coupling between two 
copies $\xi^t, \acute\xi^t$ of the Markov chain we consider here. 
To each particle we attach an additional parameter $s\in\{0,1\}$ 
called {\em state} and if a particle is paired we set its state 
to $1$ and to $0$ otherwise. We introduce also a parameter $L>0$ 
which will be used in pairing.

At $t=0$ all particles in the configurations $\xi^0$ and
$\acute\xi^0$ are at state $0$ (unpaired). At time $t>0$ first
we enumerate all particles in the configurations $\xi^t,
\acute\xi^t$ with respect to their distance to the origin. If
several particles have the same distance we enumerate them at
random (similarly later when choosing the closest particles we
shall take the one with the smallest index if there are more than
one of them). Then we start the procedure with the particle $\eta$
from the first configuration having the smallest index.

If the particle $\eta$ is at state 1 (i.e. paired) we check its
distance $\ell$ to another member of the pair $\eta'$, and if
$\ell>L$ we change the states of both these particles to $0$
(i.e. they become unpaired). Otherwise we check the distance
$\ell'$ to the closest particle $\eta''$ with state $0$ from
another process. If $\ell'<\ell$ we swap the pairing setting the
state of the particle $\eta'$ to $0$ and the state of the
particle $\eta''$ to $1$ pairing it with $\eta$.

If the particle $\eta$ is at state $0$ (may be after the
application of the previous step of the procedure) we check the
distance $\ell$ to the closest particle in the second process
$\eta'$ and if $\ell>L$ we stop the consideration of the particle
$\eta$. Otherwise if the state of $\eta'$ is $0$ we pair these
particles setting their types to $1$. If the state of $\eta'$ is
$1$ we compare $\ell$ to the distance to another member of the
paired particle $\eta'$ and proceed exactly as in the case above.

Then we continue the procedure with the next particle from the
process $\xi^t$ until all the particles will be taken into
consideration. Since only a finite number of particles may belong
to the $L$-neighborhood of a given particle the procedure is well
defined.

During the time when two particles are paired all choices of
their velocities in the coupled process are assumed to be
identical. Therefore the particles from the same pair move
synchronously until either the admissibility condition breaks down
for only one of the particles (which basically means that its
movement is blocked by another particle) or an unpaired particle
comes close enough to one of the members of the pair
(see Fig.~\ref{f:pairing}).

\Bfig(200,150)
      {\bline(0,0)(1,0)(200)   \bline(0,0)(0,1)(150)
       \bline(0,150)(1,0)(200) \bline(200,0)(0,1)(150)
       \put(10,130){\circle*{5}} \put(20,120){\circle{5}}
           \bline(10,130)(1,-1)(10)
           \put(10,135){$1$} \put(20,125){$1'$}
           \put(10,130){\vector(1,-3){9}}
           \put(20,120){\vector(1,-3){9}}
       \put(40,105){\circle*{5}} \put(40,110){$2$}
           \put(40,105){\vector(0,-1){10}}
       \put(90,120){\circle*{5}} \put(90,125){$3$}
           \put(90,120){\vector(3,-2){30}}
       \put(130,100){\circle{5}} \put(130,105){$2'$}
           \put(130,100){\vector(-2,-1){20}}
       \bline(0,75)(1,0)(200)  \put(170,135){$t$} \put(170,60){$t+1$}
       \put(20,30){\circle*{5}} \put(30,20){\circle{5}}
           \put(20,35){$1$} \put(27,25){$1'$}
       \put(40,20){\circle*{5}} \put(40,25){$2$}
           \bline(40,20)(-1,0)(10)
       \put(120,30){\circle*{5}} \put(120,35){$3$}
       \put(110,20){\circle{5}}  \put(105,25){$2'$}
           \bline(120,30)(-1,-1)(10)
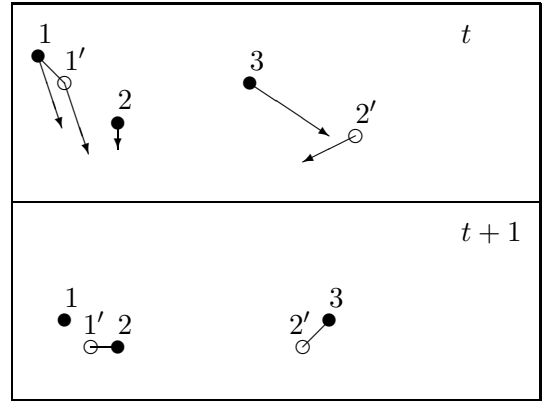
      }{Pairing of particles. Black circles corresponds to the
        particles from the first component of the process and
        open circles with primed numbers to the second component.
        The paired particles are connected by stright lines.
        At time $t$ the particles $1$ and $1'$ are paired, while
        at time $t+1$ the particle $1$ becomes unpaired while
        the particle $1'$ becomes paired with the particle $2$.
        The unpaired initially particles $3$ and $2'$ become
        paired at time $t+1$.
        The velocities at time $t$ are shown by vectors.
        \label{f:pairing}}

This construction clearly defines a Markovian coupling between
two copies of the Markov chain $\xi^t$. Of course, without
details of the procedures defining velocities and their
admissibility one cannot check the validity of the conditions
of Theorems~\ref{t:uniq-anal'} and \ref{t:uniq-prob'}.
Sufficient conditions for this will be a subject
of a forthcoming paper.

\section*{Acknowledgments}
This research has been partially supported by Russian Foundation
for Fundamental Research, CRDF and French Ministry of Education
grants.

{

}

\end{document}